\documentclass[12pt,letterpaper]{article}
\usepackage{amsfonts}
\usepackage{amssymb}
\usepackage{latexsym}  
\usepackage{pifont}  
\usepackage{amsmath}   
\usepackage{graphicx}   

\DeclareMathAlphabet{\mathbbold}{U}{bbold}{m}{n}

\let\SavedRightarrow=\Rightarrow
\usepackage{marvosym}
\let\Rightarrow=\SavedRightarrow

\newcommand{\nameX}{\raisebox{0pt}[12pt]{$\overset{\kern1pt\raisebox{-3pt}{\scriptsize $\diamond$}}{X}$}}

\newcommand\FF{{\mathcal F}}
\newcommand\BB{{\mathcal B}}
\newcommand\DD{{\mathcal D}}
\newcommand\MM{{\mathcal M}}

\newcommand\RRR{{\mathbb R}}
\newcommand\PPP{{\mathbb P}}
\newcommand\KKK{{\mathbb K}}

\newcommand\cccc{{\mathfrak c}}
\newcommand\PPPP{{\mathfrak P}}

\newcommand\cchi{{\raise 2 pt \hbox{$\chi$}}}

\newcommand\one{\mathbbold{1}} 

\newcommand\rest{\upharpoonright}     
\newcommand\res{\mathord {\upharpoonright}}  

\newcommand\onto{\twoheadrightarrow}  

\newcommand\inc{\mathrel{\raisebox{-2pt}{\rotatebox{90}{\mbox{\large $\mathbf{\models}$}}}}}

\newcommand\sinc{\mathrel{\mbox{\raisebox{-2pt}{\large \Lightning}}}}

\newcommand\MACTP{\mathrm{MA}_{\mathrm{CTP}}(\aleph_1)}  
\newcommand\MAWCTP{\mathrm{MA}_{\mathrm{WCTP}}(\aleph_1)}  
\newcommand\MA{\mathrm{MA}}
\newcommand\PFA{\mathrm{PFA}}

\newcommand\cov{\mathrm{cov}}  
\newcommand\dom{\mathrm{dom}}  
\newcommand\diam{\mathrm{diam}}   
\newcommand\next{\mathit{next}}   
\newcommand\lh{\mathrm{lh}}   
\newcommand\cdisc{\overline{D}}  

\newcommand\cat{^{\mathord{\frown}}}  

\newcommand\iv{^{-1}} 
\newcommand\lw{^{< \omega_1}} 
\newcommand\om{^{\omega_1}} 

\def\eop{{\Large \Coffeecup}}  

\newenvironment{itemizz}{\begin{itemize}\setlength{\itemsep}{-1mm}} %
{\end{itemize}}

\newtheorem{theorem}{Theorem}[section]
\newtheorem{definition}[theorem]{Definition}
\newtheorem{lemma}[theorem]{Lemma}
\newtheorem{corollary}[theorem]{Corollary}
\newtheorem{proposition}[theorem]{Proposition}

\newenvironment{proof}{{\bf Proof.}}{\eop\medskip}
\newenvironment{proofof}[1]{\medskip \textbf{Proof of #1.}}{\eop\medskip}

\begin{document}

\title{First Countable Continua and Proper Forcing}

\author{Joan E. Hart\footnote{University of Wisconsin, Oshkosh,
WI 54901, U.S.A.,
\ \ hartj@uwosh.edu}
\  and
Kenneth Kunen\footnote{University of Wisconsin,  Madison, WI  53706, U.S.A.,
\ \ kunen@math.wisc.edu}
\thanks{Both authors partially supported by NSF Grant DMS-0456653.}
}

\maketitle

\begin{abstract}
Assuming the Continuum Hypothesis,
there is a compact first countable connected space of weight $\aleph_1$
with no totally disconnected perfect subsets.
Each such space, however, may be destroyed by
some proper forcing order which does not add reals.
\end{abstract}

\section{Introduction} 
\label{sec-intro}
All topologies discussed in this paper are assumed to be Hausdorff.
As in \cite{HK2},

\begin{definition}
\label{def-weird}
A space $X$ is \emph{weird} iff $X$ is compact and not scattered,
and no perfect subset of $X$ is totally disconnected.
\end{definition}

A subset $P$  of $X$ is \emph{perfect} iff $P$ is closed and has
no isolated points.  
As usual, 
$\cccc$ denotes the (von Neumann)
cardinal $2^{\aleph_0}$.
Big weird spaces (of size $2^\cccc$) were produced from CH in
Fedorchuk, Ivanov, and van Mill \cite{FIVM}.
Small weird spaces (of size $\aleph_1$) were constructed 
from $\diamondsuit$ in \cite{HK2},
which proved:

\begin{theorem}
\label{thm-weird-HL}
Assuming $\diamondsuit$, there is a connected weird space which
is hereditarily separable and hereditarily Lindel\"of.
\end{theorem}

The weird spaces of \cite{HK2}, \cite{FIVM},
and the earlier Fedorchuk \cite{FED2} 
are all separable spaces of weight $\aleph_1$.
Our $\diamondsuit$ example is also first countable,
because it is compact and hereditarily Lindel\"of. 
In contrast,
the CH weird spaces of \cite{FIVM,FED2} have no convergent $\omega$-sequences.
We do not know whether CH can replace $\diamondsuit$ 
in Theorem \ref{thm-weird-HL}, 
but weakening hereditarily Lindel\"of to first countable we do get: 

\begin{theorem}
\label{thm-CH-weird}
Assuming CH, there is a separable first countable connected
weird space of weight $\aleph_1$.
\end{theorem}

This theorem cannot be proved by a \emph{classical} CH construction.
Classical CH arguments build the item of interest directly from
an enumeration in type $\omega_1$ of some natural set 
of size $\cccc$ (e.g.,  $\RRR$, $\RRR\lw$, etc.). 
The result, then, is preserved by any forcing which does not
add reals.  
These arguments include any CH proof found in Sierpi\'nski's text \cite{SIER}, 
as well as most CH proofs in the current literature,
including the constructions of
the big weird spaces of \cite{FED2,FIVM}.
In contrast, every space satisfying Theorem \ref{thm-CH-weird} 
is destroyed by some proper forcing order which
does not add reals.

Our proof of Theorem \ref{thm-CH-weird}
uses classical CH arguments to make $X$ weird,
but then, to make $X$ first countable, we
adapt the method of Gregory \cite{GRY} and Devlin and Shelah \cite{DS}.
The methods of \cite{GRY} and \cite{DS} are,
as Hellsten, Hyttinen, and  Shelah \cite{HHS} pointed out, 
essentially  the same.
We review the method in Section \ref{sec-pred},
and use it to prove Theorem \ref{thm-CH-weird} in Section
\ref{sec-weird}.
Although \cite{GRY} and \cite{DS} 
derive results from $2^{\aleph_0} < 2^{\aleph_1}$, 
for Theorem \ref{thm-CH-weird}, we need CH;
Section \ref{sec-rem} explains why.

In Section \ref{sec-tp}, we show that
each space satisfying Theorem \ref{thm-CH-weird} can be
destroyed by a proper forcing which does not add reals;
in $V[G]$, we add a point of uncountable character.
More precisely, if $X$ is a compactum in $V$, then in each generic extension
$V[G]$, we still have the same set $X$ with the natural topology
obtained by using the open sets from $V$ as a base.
If $X$ is first countable in $V$,  then
it must remain first countable in $V[G]$,
but $X$ need not be compact in $V[G]$.
We get the point of uncountable character
in the natural corresponding \emph{compact} 
space $\widetilde X$ in $V[G]$. 
This compact space 
determined by $X$ was
described by Bandlow \cite{BAND} (and later in \cite{DOW,DF,DK2}), 
and can be defined as follows:

\begin{definition}
\label{def-tilde}
If $X$ is a compactum in $V$ and $V[G]$ is a forcing extension of $V$,
then in $V[G]$ the corresponding compactum $\widetilde X$ 
is characterized by:
\begin{itemizz}
\item[1.] $X$ is dense in $\widetilde X$.
\item[2.] Every $f \in C(X, [0,1]) \cap V$
extends to an $\widetilde f \in C(\widetilde X, [0,1])$ in  $V[G]$.
\item[3.] The functions $\widetilde f$ \textup(for $f \in V$\textup)
separate the points of $\widetilde X$.  \samepage  
\end{itemizz}
In forcing, $\nameX$ denotes the $\widetilde X$ of $V[G]$,
while $\check X$ denotes the $X$ of $V[G]$.
\end{definition}

For example, if $X$ is the $[0,1]$ of $V$, then 
$\widetilde X$ will be the unit interval of $V[G]$;
note that in statement (2), asserted in $V[G]$, the ``$[0,1]$''
really refers to the unit interval of $V[G]$.
If in $V$, we have $X \subseteq [0,1]^\kappa$, then 
$\widetilde X$ is simply the closure of $X$ in the $[0,1]^\kappa$
of $V[G]$.
If in $V$, $X$ is the Stone space of a boolean algebra $\BB$,
then $\widetilde X$ will be the Stone space,
computed in $V[G]$, of the same $\BB$.
In general, the weights of $X$ and $\widetilde X$ will be 
the same (assuming that cardinals are not collapsed),
but their characters need not be.

Following Eisworth and Roitman \cite{ER,EIS}, we call 
a partial order $\PPP$
\emph{totally proper} iff $\PPP$ is proper and forcing with it does not
add reals.

\begin{theorem}
\label{thm-destroy}
If $X$ is compact, connected, and infinite,
and $X$ does not have a Cantor subset, then for some totally proper
$\PPP$:
$\one_{\PPP} \Vdash `` \nameX \mbox{ is not first countable''}$.
\end{theorem}

The proof is in Section \ref{sec-tp}.
Observe the importance of connectivity here.
Suppose in $V$ that $X$ is the double arrow space,
obtained from $[0,1]$ by doubling the points of $(0,1)$.
Then in any $V[G]$, $\widetilde X$ is the compactum 
obtained from $[0,1]$ by doubling the points of $(0,1) \cap V$,
and is hence first countable.

\section{Predictors} 
\label{sec-pred}

In the following,
$\lambda^{\omega_\alpha}$ denotes the set of
functions from $\omega_\alpha$ into $\lambda$.
Something like the next definition and theorem is implicit in both of
\cite{GRY,DS}:

\begin{definition}
\label{def-pred}
Let $\kappa, \lambda$ be any cardinals and $\Psi: \kappa\lw \to \lambda$.
If $f \in \kappa\om$, $g \in \lambda\om$, and  $C \subseteq \omega_1$,
then $\Psi,f$ \emph{predict} $g$ \emph{on} $C$ iff
$g(\xi) = \Psi(f \res \xi)$ for all $\xi \in C$.
$\Psi$ is a $(\kappa, \lambda)$--\emph{predictor} iff for all
$g \in \lambda\om$ there is an
$f \in \kappa\om$ and a club $C$  such that $\Psi,f$ predict $g$ on $C$.
\end{definition}

\begin{theorem}
\label{thm-pred-equiv}
The following are equivalent whenever $2 \le \kappa\le \cccc$ and
$2 \le \lambda \le \cccc$ :
\begin{itemizz}
\item[1.] There is a $(\kappa,\lambda)$--predictor.
\item[2.] There is a $(\cccc, \cccc)$--predictor.
\item[3.] $2^{\aleph_0} = 2^{\aleph_1}$.
\end{itemizz}
\end{theorem}
\begin{proof}
$(3) \to (1)$: Let $C = \omega_1 \setminus \omega$.
List $\lambda\om$ as $\{g_\alpha : \alpha < \cccc\}$, 
and choose $f_\alpha \in \kappa\om$ so that the
$f_\alpha \res \omega$, for $\alpha < \cccc$, are all distinct.
Then we can define $\Psi: \kappa\lw \to \lambda$ so that
$ \Psi(f_\alpha \res \xi) = g_\alpha(\xi)$
for all $\xi \in C$.

$(1) \to (2)$: Fix a $(\kappa,\lambda)$--predictor
$\Psi: \kappa\lw \to \lambda$.
We shall define $\Phi :  (\kappa^\omega) \lw \to (\lambda^\omega) $ so that
it is a $(\kappa^\omega,\lambda^\omega)$--predictor in the
sense of Definition \ref{def-pred}.
For $p \in (\kappa^\omega)^\xi$ and $n \in \omega$, 
define $p_{(n)} \in \kappa^\xi$ by: $p_{(n)}(\mu) = (p(\mu))(n) \in \kappa$.
Then, for $p \in (\kappa^\omega)\lw$,
define $\Phi(p) =
\left\langle \Psi(p_{(n)}) : n \in \omega \right\rangle \in \lambda^\omega$.

$(2) \to (3)$:  Fix a $(\cccc,\cccc)$--predictor $\Psi: \cccc\lw \to \cccc$.
Let $\Gamma : \cccc\lw \times \cccc\lw \to \cccc$ be any 1-1 function.
If $K \subseteq \omega_1$ is unbounded and $\xi < \omega_1$,
let $\next(\xi,K)$ be the least element of $K$ which is greater than $\xi$.

For each $B \in \cccc\om$, choose $G(n,B), F(n,B) \in \cccc\om$ and
clubs $C(n,B) \subseteq \omega_1$ for $n \in \omega$ as follows:
Let $G(0,B) = B$.
Given $G(n,B)$, let $C(n,B)$ be club of limit ordinals
and let $F(n,B) \in \cccc\om$ be such that
$(G(n,B) )(\xi) = \Psi((F(n,B) ) \rest \xi)$ for all
$\xi \in C(n,B)$.
Then define $G(n+1,B)$ so that
\[
(G(n+1,B) )(\xi) = \Gamma\big(
F(n,B) \rest \next(\xi, C(n,B)) ,
\ G(n,B) \rest \next(\xi, C(n,B))  \big) \ \ 
\]
for each $\xi$.

Now, fix $B,B' \in \cccc\om$, and consider the statement:
\[
\forall n\in \omega  \,  \big[ G(n,B) \res \xi = G(n,B') \res \xi \big]  
\tag{\ding{"50}$(\xi)$}
\]
So, \ding{"50}$(0)$ is true trivially, and 
\ding{"50}$(\xi)$ implies 
\ding{"50}$(\zeta)$ whenever $\zeta < \xi$.
We shall prove inductively that
\ding{"50}$(1)$ implies 
\ding{"50}$(\eta)$ for all $\eta < \omega_1$.
If we do this, then \ding{"50}$(1)$ will imply $B = B'$, so we shall
have $2^{\aleph_0} = 2^{\aleph_1}$, since there are
$2^{\aleph_1}$ possible values for $B$ but only 
$2^{\aleph_0}$ possible values for the sequence
$\langle (G(n,B))(0) : n \in \omega \rangle$.

The induction is trivial at limits, so it is sufficient to
fix $\eta$ with $1 \le \eta < \omega_1$,
assume \ding{"50}$(\eta)$, and prove
\ding{"50}$(\eta + 1)$ --- that is,
$(G(n,B))(\eta)  = (G(n,B'))(\eta)$ for all $n$.
Fix $n$.  
For $\xi < \eta$,
we have 
$(G(n+1,B))(\xi)  = (G(n+1,B'))(\xi)$, which implies:
\begin{itemizz}
\item[a.] $\next(\xi, C(n,B))  = \next(\xi, C(n,B'))$;
call this $\gamma_\xi$.
\item[b.] $F(n,B) \res \gamma_\xi = F(n,B') \res \gamma_\xi$.
\item[c.] $G(n,B) \res \gamma_\xi = G(n,B') \res \gamma_\xi$.
\end{itemizz}
Applying (a) for all $\xi < \eta$: $\eta \in C(n,B)$ iff $\eta \in C(n,B')$.
If $\eta \notin C(n,B),C(n,B')$, then fix $\xi$ with
with $\xi < \eta < \gamma_\xi$; now (c) implies
$(G(n,B))(\eta)  = (G(n,B'))(\eta)$.
If $\eta \in C(n,B),C(n,B')$, then $\eta$ is a limit ordinal
and (b) implies $F(n,B) \res \eta = F(n,B') \res \eta$;
now $(G(n,B))(\eta)  = (G(n,B'))(\eta) = \Psi((F(n,B) ) \rest \eta)$.
\end{proof}

The \emph{non}-existence of a $(2,2)$--predictor 
is the weak version of $\diamondsuit$ discussed by
Devlin and Shelah in \cite{DS}, where they use it to prove
that, assuming $2^{\aleph_0} < 2^{\aleph_1}$,
every ladder system on $\omega_1$ has a non-uniformizable coloring.
By Shelah \cite{SH} (p.~196),
each such coloring may be uniformized
in some totally proper forcing extension.

A direct proof of $(3) \to (2)$, resembling the above proof of
$(3) \to (1)$, would obtain $C$ fixed at
$\omega_1 \setminus \{0\}$, since one may choose the $f_\alpha$
so that the $f_\alpha(0)$, for $\alpha < \cccc$, are all distinct.
Gregory \cite{GRY} used the failure of $(2)$, with this specific $C$,
to derive a result about trees under $2^{\aleph_0} < 2^{\aleph_1}$;
see Theorem \ref{thm-gregory} below.

\section{Some Totally Proper Orders} 
\label{sec-tp}
We consider \emph{forcing posets}, $(\PPP ; \le, \one)$,
where $\le$ is a transitive and reflexive relation on $\PPP$
and $\one$ is a largest element of $\PPP$.
As usual, if $p,q \in \PPP$,  then
$p \not\perp q$ means that $p,q$ are compatible (that is, have a common
extension),
and $p \perp q$ means that $p,q$ are incompatible.

\begin{definition}
Assume that $X$ is compact, connected, and infinite.
Let $\KKK = \KKK_X$ be the forcing poset consisting
of all closed, connected, infinite subsets of $X$,
with $p \le q$ iff $p \subseteq q$ and $\one_\KKK = X$.
In $\KKK$, define $p \inc q$ iff $p \cap q = \emptyset$.
\end{definition}

Note that $p \perp q$ iff $p \cap q$ is totally disconnected.
The stronger relation $p \inc q$ will be useful in the proof that
$\KKK$ is totally proper whenever $X$ does not have a Cantor subset.
First, we verify that $\KKK$ is separative; this follows easily
from the following lemma, which is probably well-known; a proof is
in \cite{HK2}:

\begin{lemma}
\label{lemma-con-split}
If $P$ is compact, connected, and infinite, and  $U \subseteq P$
is a nonempty open set, then there is a closed $R \subseteq U$
such that $R$ is connected and infinite.
\end{lemma}

In particular, in $\KKK$, if $p \not\le q$, then we may apply
this lemma with $U = p \setminus q$ to get $r \le p$
with $r \perp q$, proving the following: 

\begin{corollary}
\label{cor-sep}
If $X$ is compact, connected, and infinite, then $\KKK_X$ is
separative and atomless.
\end{corollary}

We collect some useful properties of
the relation $\inc$ on $\KKK$ 
in the following:

\begin{definition}
A binary relation $\sinc$ on a forcing poset is a
\emph{strong incompatibility relation} iff
\begin{itemizz}
\item[1.]
$p \sinc q$ implies $p \perp q$.
\item[2.]
Whenever $p \perp q$, there are $p_1,q_1$ with
$p_1 \le p$,  $q_1 \le q$, and $p_1 \sinc q_1$.
\item[3.]
$p \sinc q \ \&\ p_1 \le p \ \&\ q_1 \le q \ \to\ p_1 \sinc q_1$.
\end{itemizz}
\end{definition}

This definition does not require $\sinc$ to be symmetric,
but note that the relation $p \sinc q \ \&\ q \sinc p$
is symmetric and is also a strong incompatibility relation.

\begin{lemma}
The relation $\inc$ is a strong incompatibility relation on $\KKK_X$.
\end{lemma}
\begin{proof}
Conditions (1) and (3) are obvious.  For (2):
Suppose that $p \perp q$.  Let $F = p \cap q$, which is totally
disconnected.
Then by Lemma \ref{lemma-con-split}
there is an infinite
connected $p_1 \subseteq p \backslash F$.  
Likewise, we get
$q_1 \subseteq q \backslash F$.  
\end{proof}

\begin{definition}
If $\PPP$ is a forcing poset with a strong incompatibility 
relation $\sinc$, then a \emph{strong Cantor tree} in $\PPP$
\textup(with respect to $\sinc$\textup) is a subset
$\{p_s : s \in 2^{<\omega}\} \subseteq \PPP$ such
that each $p_{s\cat\mu} < p_s$ for $\mu = 0,1$,
and each $p_{s\cat 0}\sinc p_{s\cat 1}$.
Then, $\PPP$ has the \emph{weak Cantor tree property (WCTP)}
\textup(with respect to $\sinc$\textup)
iff whenever $\{p_s : s \in 2^{<\omega}\} \subseteq \PPP$
is a strong Cantor tree, there is at least one $f \in 2^\omega$ such that
$\PPP$ contains some $q = q_f$ with
$q \le p_{f\res n}$ for each $n \in\omega$.
\end{definition}

Note that if $\PPP$ has the 
WCTP, then the set of $f$ for which $q_f$ is defined must
meet every perfect subset
of the Cantor set $2^\omega$, since otherwise we could find a subtree
of the given Cantor tree which contradicts the WCTP.

\begin{lemma}
\label{lemma-compact-WCTP}
If $X$ is compact, connected, and infinite,
and $X$ does not have a Cantor subset,
then $\KKK_X$ has the WCTP.
\end{lemma}

\begin{definition}
$\PPP$ has the \emph{Cantor tree property (CTP)} iff
$\PPP$ has the WCTP with respect to the usual $\perp$ relation.
\end{definition}

$\KKK_X$ need not have the CTP (see Theorem \ref{thm-CH-weirder}).
A countably closed $\PPP$ clearly has the CTP.
In the case of trees, the
CTP was also discussed in \cite{HHS}
(where it was called ``$\aleph_0$ fan closed'') and in \cite{HK2}.
The following modifies Lemma 3 of \cite{HHS} and Lemma 5.5 of \cite{HK2}:

\begin{lemma}
\label{lemma-tot-prop}
If  $\PPP$ has the WCTP, then $\PPP$ is totally proper.
\end{lemma}
\begin{proof} 
Define $q \le' p$ iff there is no $r$ such that $r \le q$ and
$r \perp p$.
When $\PPP$ is separative, this is equivalent to $q \le p$.

Fix a suitably large regular cardinal $\theta$, and let $M \prec H(\theta)$
be countable with $(\PPP ; \le, \one, \sinc) \in M$,
and fix $p \in \PPP \cap M$.  It suffices (see \cite{ER})
to find a $q \le p$ such that whenever $A \subseteq \PPP$
is a maximal antichain and $A \in M$, there is an $r \in A\cap M$ with
$q \le' r$.
If $\PPP$ has an atom $q \le p$ such that $q \in M$,
then we are done.  Otherwise, then
since $M \prec H(\theta)$, $\PPP$ must be atomless below $p$.
Let $\{A_n : n \in \omega\}$ list all the maximal
antichains which are in $M$.  Build a strong Cantor tree
$\{p_s : s \in 2^{<\omega}\} \subseteq \PPP \cap M$ such that,
$p_{()} \le p$, and such that,
when $n\in \omega$ and $s \in 2^n$,
$p_s$ extends some element of $A_n \cap M$. 
Then choose $f \in 2^\omega$ such that there is some $q \in \PPP$ with
$q \le p_{f\res n}$ for each $n \in\omega$.
\end{proof}

\begin{proofof}{Theorem \ref{thm-destroy}}
Let $\PPP = \KKK_X$.
Working in $V[G]$, let $G' = \{\widetilde p : p \in G\}$;
then $\bigcap G' = \{y\}$ for some  $y \in \widetilde X\setminus X$.
Since $\PPP$ does not add $\omega$--sequences,
$\bigcap E \supsetneqq \{y\}$ whenever $E$ is a countable subset
of $G'$.  Thus, $\cchi(y, \widetilde X)$ is uncountable.
\end{proofof}

These totally proper partial orders yield
natural weakenings of PFA:

\begin{definition}
\label{def-MACTP}
If $\PPPP$ is a class of forcing posets, then
$\MA_\PPPP(\aleph_1)$ is the statement that whenever $\PPP \in \PPPP$
and $\DD$ is a family of $\le \aleph_1$ dense subsets of $\PPP$,
then there is a filter on $\PPP$ meeting each $D \in \DD$.
\end{definition}

Trivially, $\PFA \to \MAWCTP \to \MACTP$,
but in fact $\MAWCTP \leftrightarrow \MACTP$
(see Lemma \ref{lemma-mactp-equiv}).
Also, $\MACTP \to 2^{\aleph_0} = 2^{\aleph_1}$ 
(see Corollary \ref{cor-gregory}),
so, the natural iteration of (totally proper)
CTP orders with countable supports must introduce reals
at limit stages.
By the proof of Theorem 5.9 in \cite{HK2},
PFA does not follow from
$\MACTP + \MA(\aleph_1) + 2^{\aleph_0} = \aleph_2$,
which in fact can be obtained by ccc forcing over $L$.

We now consider some CTP trees.

\begin{definition}
\label{def-order}
Order $\lambda\lw$ by: $p \le q$ iff $p \supseteq q$.
Let $\one = \emptyset$, the empty sequence.
\end{definition}

So, $\lambda\lw$ is a tree, with the root $\one$ at the top.
Viewed as a forcing order, it is equivalent to countable partial
functions from $\omega_1$ to $\lambda$. 
We often view $p \in \lambda\lw$ as a countable sequence and
let $\lh(p) = \dom(p)$.  Then $\lh(\one) = 0$.

Kurepa showed that SH is equivalent to the non-existence of
Suslin trees.  A similar proof shows that $\MACTP$ is
equivalent to the non-existence of Gregory trees:

\begin{definition}
A \emph{Gregory tree} is a forcing poset $\PPP$
which is a subtree of $\cccc\lw$ and satisfies:
\begin{itemizz}
\item[1.] $\PPP$ has the CTP.
\item[2.] $\PPP$ is atomless.
\item[3.] $\PPP$ has no uncountable chains.
\end{itemizz}
\end{definition}

It is easily seen that if any of conditions (1)(2)(3) are dropped,
such trees may be constructed in ZFC.  However:

\begin{lemma}
\label{lemma-mactp-equiv}
The following are equivalent:
\begin{itemizz}
\item[1.] $\MACTP$.
\item[2.] $\MAWCTP$.
\item[3.] There are no Gregory trees.
\end{itemizz}
\end{lemma}
\begin{proof}
$(1) \rightarrow (3)$:  Let $\PPP$ be a Gregory tree.
As with Suslin trees under $\MA(\aleph_1)$,  a filter $G$
meeting the sets $D_\xi := \{p \in \PPP : \lh(p) \ge \xi\}$
yields an uncountable chain, and hence a contradiction,
but to apply $\MACTP$, we must prove that each $D_\xi$ is dense in $\PPP$.
To do this, induct on $\xi$.  The case $\xi = 0$ is trivial.
For the successor stages, use the fact that $\PPP$ is atomless.
For the limit stages, use the CTP.

$(3) \rightarrow (2)$:  Fix $\PPP$ with the WCTP and dense
sets $D_\xi \subseteq \PPP$  for $\xi < \omega_1$.
We need to produce a filter $G \subseteq \PPP$ meeting each $D_\xi$.
This is trivial if $\PPP$ has an atom, so assume that $\PPP$
is atomless.

Inductively define a subtree $T$ of $2\lw$ together with a
function $F : T \to \PPP$ as follows:
$F(\one) = \one_\PPP$.  If $t\in T$ and $\lh(t) = \xi$,
then $t\cat 0 \in T$ and $t \cat 1 \in T$,
and $F(t\cat 0)$, $F(t\cat 1)$ are extensions of 
$F(t)$ such that each $F(t\cat i) \in D_\xi$
and $F(t\cat 0) \sinc F(t\cat 1)$;
to accomplish this, given $t$ and $F(t)$: first choose two
$\perp$ extensions of $F(t)$, then extend these to be $\sinc$,
and then extend these to be in $D_\xi$.
If $\eta < \omega_1$ is a limit ordinal and $\lh(t) = \eta$,
then $t \in T$ iff
$\forall \xi < \eta \, [ t \res \xi \in T]$ and
$\exists q \in \PPP \, \forall \xi < \eta \, [ q \le F( t\res \xi) ]$;
then choose $F(t)$ to be some such $q$.

$T$ is clearly atomless, and $T$ has the CTP because $\PPP$
has the WCTP.  If there are no Gregory trees, then $T$ has
an uncountable chain, so fix $g \in 2\om$ such that
$g \res \xi \in T$ for all $\xi < \omega_1$, and
let $G = \{y \in \PPP : \exists \xi < \omega_1\, [F(g\res \xi) \le y]\}$.
\end{proof}

\begin{theorem}[Gregory \cite{GRY}]
\label{thm-gregory}
If  $2^{\aleph_0} < 2^{\aleph_1}$ then there is a Gregory tree.
\end{theorem}

\begin{corollary}
\label{cor-gregory}
$\MACTP$ implies that $2^{\aleph_0} = 2^{\aleph_1}$.
\end{corollary}

\section{A Weird Space}
\label{sec-weird}

We now prove Theorem \ref{thm-CH-weird}.
The basic construction is an inverse limit in $\omega_1$ steps,
and we follow approximately the terminology in
\cite{DK1, HK2}.  We build a compact space
$X_{\omega_1} \subseteq [0,1]\om$
by constructing inductively
$X_\alpha \subseteq [0,1]^{1 + \alpha} \cong [0,1]\times [0,1]^\alpha$.
Usually, one has $X_\alpha \subseteq [0,1]^{\alpha}$ in these
constructions, but for finite $\alpha$, the notation will be slightly
simpler if we start at stage $0$ with $X_0 = [0,1] = [0,1]^1$;
of course, $1 + \alpha = \alpha$ for infinite $\alpha$.

\begin{definition}
$\pi^\beta_\alpha: [0,1]^{1 +\beta} \onto [0,1]^{1 +\alpha}$
is the natural projection. 
\end{definition}

As usual, $\pi : X \onto Y$ means that $\pi$
is a \emph{continuous} map from $X$ \emph{onto} $Y$.
These constructions
always have $\pi^{ \beta}_{ \alpha} (X_\beta) =  X_\alpha$
whenever $0 \le \alpha \le \beta \le \omega_1$.
This determines $X_\gamma$ for limit
$\gamma$, so the meat of the construction involves describing
how to build $X_{\alpha+1}$ given $X_\alpha$.

A classical CH argument can ensure that $X_{\omega_1}$ is
weird, but by Theorem \ref{thm-destroy}, such an argument cannot make
$X_{\omega_1}$ first countable.
However, the same classical argument will let us construct a binary tree
of spaces, resulting in a weird space $X_g \subseteq [0,1]\om$
for each $g \in 2\om$.  We shall show that if no
$X_g$ were first countable, then there would be
a $(\cccc,2)$--predictor $\Psi: [0,1]\lw \to 2$; 
so CH ensures that some $X_g$ is first countable.

Our tree will give us an $X_p$ for each $p \in 2^{\le \omega_1}$. 
We now list requirements (R1)(R2)(R3)$\cdots$(R17) on the construction;
a proof that all the requirements can be satisfied,  
and that they yield a weird space, concludes
this section.
We begin with the requirements involving the inverse limit:

\begin{itemizz}
\item[R1.] $X_\one = [0,1]$, where $\one$ is the empty sequence.
\item[R2.] $X_p$ is an infinite closed connected
subspace of $[0,1]^{1 + \lh(p)}$.
\item[R3.] $\pi^{ \beta}_{ \alpha} \res X_p: X_p \onto X_{p \res \alpha}$, and
is irreducible, whenever $\beta = \lh(p) \ge \alpha$.
\end{itemizz}
When $\gamma  = \lh(p) \le \omega_1$ is a limit, (R2)(R3) force:
\[
X_p = 
\{x \in [0,1]^\gamma : \forall \alpha < \gamma \, 
[\pi^{\gamma}_{\alpha}(x) \in X_{p \res \alpha} ]\} \ \ .
\tag{\ding{"60}}
\]
To simplify notation for the restricted projection maps, we shall use:

\begin{definition}
If $\beta = \lh(p) \ge \alpha$ and $r = p \res \alpha$, define
$\pi^p_r = \pi^\beta_\alpha \res X_p : X_p \onto X_r$.
\end{definition}

As in \cite{HK2}, each of $X_{p\cat 0}$ and $X_{p\cat 1}$
is obtained from $X_p$ as the graph of a ``$\sin(1/x)$'' curve.
We choose $h_q, u_q$, and $v_q^n$
for $n < \omega$ and $q \in 2\lw$ of successor length, satisfying,
for $i = 0,1$:

\begin{itemizz}
\item[R4.] $u_{p\cat i} \in X_{p}$ and
$h_{p\cat i} \in C(X_p \backslash \{u_{p \cat i}\},\; [0,1])$
and $X_{p \cat i} = \overline{h_{p\cat i}}$.
\item[R5.] $v_{p\cat i}^n \in  X_{p} \backslash\{u_{p\cat i}\}$, and
$\langle v_{p\cat i}^n: n \in \omega \rangle \to u_{p\cat i}$, and all points
of $[0,1]$ are limit points of
$\langle h_{p\cat i}(v_{p\cat i}^n): n \in \omega \rangle$.
\end{itemizz}
As usual, we identify $h_{p\cat i}$ with its graph.
So, if $\alpha = \lh(p)$, then $X_{p \cat i}$ is a subset
of  $[0,1]^{1 + \alpha} \times [0,1]$, which we identify with
$[0,1]^{1 + \alpha + 1}$.
We shall say that the point $u_{p\cat i}$ gets
\emph{expanded} in the passage from $X_p$ to $X_{p\cat i}$;
the other points get \emph{fixed}.
(R3) follows from (R4) plus (\ding{"60}).
Also, if $\delta < \alpha$, then
$\pi^{p}_{p \res \delta}: X_p \onto X_{p \res \delta}$,
and $(\pi^{ p}_{p \res \delta} )\iv\{x\}$ is a singleton unless $x$ is in
the countable set
$\{\pi^{p \res \xi}_{p \res\delta}(u_{p\res(\xi + 1})) :
\delta \le \xi < \alpha\}$.

We now explain how points in $X_g \subset [0,1]\om$ 
can predict $g$, in the sense of Definition \ref{def-pred}.
We shall get $A_q$ and $B_q$ for $q \in 2\lw$ of successor length,
satisfying:

\begin{itemizz}
\item[R6.]  For $i = 0,1$: $A_{p\cat i},  B_{p\cat i} \subseteq X_p$ and
$A_{p\cat i} = X_p \setminus B_{p\cat i}$.  \samepage  
\item[R7.]  For $i = 0,1$ and $\xi < \lh(p)$ :
$A_{p\cat i} \supseteq (\pi^p_{p \res\xi})\iv (A_{p \res (\xi + 1)}) $.
\item[R8.] $ B_{p\cat 0} \cap B_{p\cat 1} = \emptyset$
\item[R9.]  For $i = 0,1$: $u_{p \cat i} \in B_{p\cat i}$.
\end{itemizz}

Observe that some care must be exercised here in the inductive
construction; otherwise, at some stage (R7) might imply that
$A_{p\cat i} =  X_p$, so that 
$B_{p\cat i} =  \emptyset$, making (R9) impossible.

(R6)(R7)(R9) imply that points in $ A_{p\cat i}$ are forever fixed
in the passage from $X_p$ to any future $X_q$ with $q \le p\cat i$;
only points in $B_{p\cat i}$ can get expanded.
Points which are forever fixed must wind up having countable character,
and (R8) lets us use a point of uncountable character in $X_g$ to
predict $g$:

\begin{lemma}
\label{lemma-first-count}
Assume that we have {\rm (R1 -- R9)}, and assume
that $2^{\aleph_0} < 2^{\aleph_1}$.
Then $X_g$ is first countable for some $g \in 2\om$.
\end{lemma}
\begin{proof}
We shall define $\Psi: [0,1]\lw \to 2$, and prove
that $\Psi$ is a $(\cccc,2)$--predictor if every $X_g$ contains a point of
uncountable character.

Say $\lh(p) = \alpha < \omega_1$ and $\delta < \alpha$.
If $x \in B_{p\cat i} \subseteq X_p$ then, by (R6)(R7),
$\pi^\alpha_\delta(x) \in B_{p\res(\delta + 1)} \subseteq X_{p \res \delta}$.
Applying (R8),
if $x \in [0,1]^{1 + \alpha}$ and 
$x \in B_{p\cat i} \cap B_{r\cat j}$, then $p = r$ and $i = j$;
to prove this, consider the least $\delta < \alpha$
such that $p(\delta) \ne r(\delta)$.

Set $\Psi(x) = 0$ if $\lh(x) < \omega$.
Now, say $x \in [0,1]^\alpha$, where $\omega \le \alpha < \omega_1$
(so $1 + \alpha = \alpha$).
If there exist $p \in 2^\alpha$ and $i \in 2$ such
that $x \in B_{p\cat i}$, then these $p,i$ are unique, and set $\Psi(x) = i$.
If there are no such $p,i$, then set $\Psi(x) = 0$.

Now, assume that for each $g$, we can find $z = z_g \in X_g$
with $\cchi(z, X_g) = \aleph_1$.
Let $C = \omega_1 \setminus \omega$.
We shall show that $\Psi,z$ predict $g$ on $C$.
For $\xi \in C$, let $p\cat i = g \res (\xi + 1) $.
Then $z \res \xi = \pi^g_p(z)  \in X_p$, and
$z \res \xi$ must be in $B_{p\cat i}$, since if it were in $A_{p\cat i}$, then
$(\pi^g_p)\iv(\pi^g_p(z)) = \{z\}$, so that $\cchi(z, X_g) = \aleph_0$.
Thus, $\Psi(z \res \xi) = i = g(\xi)$.
\end{proof}

Since every $X_g$ clearly has weight $\aleph_1$, we are done
if we can make every $X_g$ weird.
Since points in $A_{p\cat i}$ are forever fixed,
we must make sure that $A_{p\cat i}$ has no Cantor subsets.
Conditions (R6)(R8) say that
$ A_{p\cat 0} \cup A_{p\cat 1} = X_p$, so 
$ A_{p\cat 0}$ and  $A_{p\cat 1}$ must be Bernstein sets.
Note that Condition (R7) may present a problem at limit stages.
When $\lh(p) = \alpha$  we have
$A_{p\cat i} \supseteq
\bigcup_{\xi < \alpha}(\pi^p_{p \res \xi})\iv (A_{p \res (\xi + 1)}) $.
Points in $A_{p \res (\xi + 1)}$ are forever fixed,
so each $(\pi^p_{p \res \xi})\iv (A_{p \res (\xi + 1)})$ will have no
Cantor subsets. 
Without further requirements, though, 
$\bigcup_{\xi < \alpha}(\pi^p_{p \res \xi})\iv (A_{p \res (\xi + 1)}) $
may contain a Cantor subset.
So, we make sure each such union is disjoint from some set 
in a tree of Bernstein sets:

\begin{definition}
\label{def-bern}
For any topological space $Y$ and $p \in 2\lw$,
a \emph{Bernstein tree in} $Y$ \emph{rooted in} $p$ is a family
of subsets of $Y$, $\{D^q : q \le p\}$, satisfying:
\begin{itemizz}
\item[1.] For each $q$, neither 
$D^q$ nor $Y \backslash D^q$ contains a Cantor subset.
\item[2.] Each $D^{q\cat 0}  \cap D^{q\cat 1} = \emptyset$.
\item[3.] If $r \le q$ then $D^r \subseteq D^q$.
\end{itemizz}
\end{definition}

Note that if $Y$ itself does not contain a Cantor subset, then (1)
is trivial, and we may take all $D^q = \emptyset$ to satisfy (2) and (3).

Now, in our construction, we also
build $D_p^q$ for $q \le p \in 2\lw$ satisfying:

\begin{itemizz}
\item[R10.] For each $p \in 2\lw$: $\{D_p^q : q \le p\}$ is
a Bernstein tree in $X_p$ rooted in $p$.
\item[R11.] If $q \le p \le r$ and
$\pi = \pi^{p}_{r} : X_p \onto X_r$ and
$x \in X_p$ with $\pi\iv(\pi(x)) = \{x\}$, then
$x \in D_p^q$ iff $\pi(x) \in D_r^q$.
\item[R12.] For each $p \in 2\lw$ and $i \in 2$:
$B_{p\cat i} = D_p^{p\cat i}$ and
$A_{p\cat i} = X_p \setminus D_p^{p\cat i}$.
\end{itemizz}

Of course, (R12) simply defines $A_{p\cat i}$ in terms of the
$D_p^q$, and then (R10) guarantees that no $A_{p\cat i}$ has a
Cantor subset, but we need to verify that the conditions
(R1 -- R12) can indeed be satisfied.  First, three easy lemmas
about Bernstein trees.
A standard inductive construction in $\cccc$ steps shows:

\begin{lemma}
\label{lemma-bern-exists}
If $Y$ is a separable metric space, then
there is a Bernstein tree in $Y$ rooted in $\one$.
\end{lemma}

Using the fact that every uncountable Borel subset of the 
Cantor set contains a perfect subset, we get:

\begin{lemma}
\label{lemma-split-bern}
Assume that $Y$ is any topological space,
$Z$ is a Borel subset of $Y$,
and $\{D^q : q \le p\}$ is a family of subsets of $Y$
satisfying $(2)(3)$ of Definition \ref{def-bern}.
Then 
$\{D^q : q \le p\}$ is a Bernstein tree in $Y$ iff both
$\{D^q \cap Z : q \le p\}$ is a Bernstein tree in $Z$ and
$\{D^q \backslash Z : q \le p\}$ is a Bernstein tree in $Y \backslash Z$.
\end{lemma}

Combining these two lemmas:

\begin{lemma}
\label{lemma-add-bern}
If $Y$ is a separable metric space, $Z$ is a Borel subset of $Y$,
and $\{E^q : q \le p\}$ is a Bernstein tree in $Z$ rooted in $p$,
then there is a Bernstein tree $\{D^q : q \le p\}$ 
in $Y$ rooted in $p$ such that each $D^q  \cap Z = E^q$.
\end{lemma}

Returning to the construction:

\begin{lemma}
\label{lemma-cond}
There exist $X_p$ for $p \in 2^{\le \omega_1}$
satisfying Conditions {\rm (R1 -- R12)}.
\end{lemma}
\begin{proof}
We start with $X_\one = [0,1]$, and 
we obtain the $D_\one^q$ by applying Lemma \ref{lemma-bern-exists}.

If $\alpha = \lh(p) > 0$ and
we have done the construction for $p \res \xi$
for all $\xi < \lh(p)$, then $X_p$ is determined either by 
(R4) when $\lh(p)$ is a successor or by
(\ding{"60}) when $\lh(p)$ is a limit.
If $\alpha < \omega_1$, we 
construct the $D_p^q$ to satisfy (R10)(R11) as follows:
For $\xi < \alpha$, use $\pi_\xi$ for $\pi^{p}_{p \res\xi}$.  Let
$Z_\xi  = \{x \in X_p : \pi_\xi \iv(\pi_\xi (x)) = \{x\} \}$,
and let $Z = \bigcup_{\xi < \alpha} Z_\xi$.
Observe that $Z$ and all the $Z_\xi$ are Borel sets.
Let $\{E_p^q : q \le p\}$ be the Bernstein tree in $Z$ rooted in $p$
defined by saying that for $x \in Z_\xi$:
$x \in E_p^q$ iff $\pi_\xi(x) \in D_{p\res\xi}^q$.  Note that, by (R11) applied
inductively, this is independent of which $\xi$ is used.
To obtain the $D_p^q$ from the $E_p^q$, 
apply Lemma \ref{lemma-add-bern}.
Note that, by (R11) applied inductively once again,
these $D_p^q$ work for $X_p$.

The $ A_{p\cat i}$ and $ B_{p\cat i}$ (for $i = 0,1$) are now defined by
(R12), and we must verify that this definition satisfies (R7):
Assume that  $\xi < \lh(p) = \alpha$  and $x \in X_p$ and
$\pi_\xi(x) \in A_{p \res (\xi + 1)} $.
We must show that $x \in A_{p\cat i}$; equivalently,
by (R12), that $x \notin D_p^{p\cat i}$.
Now $\pi_\xi(x) \in A_{p \res (\xi + 1)} $ implies that
$\pi_\xi\iv(\pi_\xi(x)) = \{x\}$ (using (R4 -- R9) inductively),
so that $x \notin D_p^{p\cat i}$ iff 
$\pi_\xi(x) \notin D_{p\res \xi}^{p\cat i}$.
By (R10) for $p\res\xi$
and Definition \ref{def-bern}(3),
$D_{p\res \xi}^{p\cat i} \subseteq D_{p\res \xi}^{p\res (\xi + 1)}$.
So $A_{p \res (\xi + 1)} = X_{p \res\xi} \setminus 
D_{p\res \xi}^{p\res (\xi + 1) }$ gives us (R7).

Since the $B_{p\cat i}$ are nonempty, there is no problem
choosing the $u_{p\cat i}$, $v_{p\cat i}^n$, and $h_{p\cat i}$
to satisfy (R4)(R5)(R9), and then the $X_{p\cat i}$ are 
defined by (R4).
\end{proof}

Finally, we must make each $X_g$ weird.  Observe:

\begin{lemma}
\label{lemma-conn}
Conditions {\rm (R1 -- R5)} imply that
if $F \subseteq X_p$ is closed and connected then
$ (\pi^{q}_{p}) \iv (F)$ is connected
for all $q \le p$.
\end{lemma}

Now, we shall make sure that whenever $F$ is a perfect subset of
$X_g$, there is some $\alpha < \omega_1$ such that
$ (\pi^g_{g \res(\alpha + 1)})\iv ( \{u_{g \res (\alpha + 1)} \} \times [0,1])
\subseteq F$ (recall that our construction gave us 
$\{u_{g \res (\alpha + 1)} \} \times [0,1] \subset 
X_{g \res (\alpha + 1)} \subset X_{g \res \alpha} \times [0,1] $) .
By Lemma \ref{lemma-conn}, this implies 
that $F$ is not totally disconnected.
The argument in \cite{HK2} obtained this $\alpha$ by using
$\diamondsuit$ to capture $F$.  Here, we 
replace this use of $\diamondsuit$ by a classical CH argument.
First, as in \cite{HK2},
construct $\FF_p$ for $p \in 2\lw$ so that:

\begin{itemizz}
\item[R13.] $\FF_p$ is a countable family of uncountable
closed subsets  of $X_p$.
\item[R14.] If $F \in \FF_p$ and $q \le p$ then
$ (\pi^{q}_{p})\iv(F) \in \FF_{q}$.
\item[R15.] For each $F \in \FF_p$, either $u_{p \cat i} \notin F$,
or $u_{p \cat i} \in F$ and $v_{p \cat i}^n \in F$ for all
but finitely many $n$.
\item[R16.]
$\{u_{p\cat i}\} \times [0,1] \in \FF_{p \cat i}$.
\end{itemizz}

We may satisfy (R13)(R14)(R16) simply by defining 
\[
\FF_p = 
\left\{(\pi^p_{p \res \xi })\iv \{ u_{p \res(\xi + 1)} \right\} :
\xi < \lh(p) \} \ \ .
\]
Requirements (R4)(R14)(R15) imply:

\begin{lemma}
\label{lemma-irred}
$\pi^q_p :  (\pi^{q}_{p})\iv(F) \onto F$
is irreducible whenever $F \in \FF_p$  and $q \le p$.
\end{lemma}

Then, we use CH rather than $\diamondsuit$ to get:

\begin{itemizz}
\item[R17.]  Whenever $p \in 2\lw$ and $F$ is an uncountable
closed subset of $X_p$, there is a $\beta$ with
$\lh(p) < \beta < \omega_1$ such that for all 
$q < p$ with $\lh(q)  = \beta$ 
and for each $x \in \{u_{q \cat 0}, u_{q \cat 1}\} 
\cup \{v_{q \cat i}^n : n \in \omega \ \&\ i \in 2\}$,
the projections
$\pi = \pi^{q}_{p}$ satisfy
$\pi(x) \in F$ and $|\pi\iv(\pi(x))| = 1$.
\end{itemizz}

\begin{proofof}{Theorem \ref{thm-CH-weird}}
Assuming that we can obtain (R1 -- R17), 
note that each $X_g$ is separable, because
each $\pi^g_\one : X_g \onto X_\one$ is irreducible.
Then, to finish, by Lemma \ref{lemma-first-count}, 
it suffices to show that each $X_g$ is weird.
Fix a perfect $H \subseteq X_g$;
we shall show that it is not totally disconnected.
First, fix $\alpha < \omega_1$ such that, if we set
$p = g\res \alpha$ and $F = \pi^g_p(H)$, then $F$ is perfect
(the set of all such $\alpha$ form a club).
Then, fix $\beta > \alpha$ as in (R17), let $q = g\res\beta$,
and let $i = g(\beta)$, so that $q \cat i = g\rest(\beta+1)$.
Let $K  = (\pi^q_p)\iv(F)$.
Then $\pi^g_q(H) \subseteq K$, and this inclusion may well be proper.
However, 
$u_{q \cat i} \in \pi^g_q(H)$ and $v_{q \cat i}^n \in \pi^g_q(H)$
for each $n \in \omega$ because
$\pi(u_{q \cat i}) \in F$ and $\pi(v_{q \cat i}^n) \in F$
and $|\pi\iv(\pi(u_{q \cat i}))|  = |\pi\iv(\pi(v_{q \cat i}^n))| = 1$.
It follows (using (R5)) that
$E := \{u_{q\cat i}\} \times [0,1] \subseteq \pi^g_{q\cat i}(H)$.
Since $E \in \FF_{q \cat i}$  by (R16) and $H$ maps onto $E$,
Lemma \ref{lemma-irred} implies that $(\pi^{g}_{q\cat i})\iv(E) \subseteq H$.
Since $(\pi^{g}_{q\cat i})\iv(E)$ is connected by Lemma
\ref{lemma-conn}, $H$ cannot be totally disconnected.

Next, to obtain conditions (R1 -- R17),
we must augment the proof of Lemma \ref{lemma-cond}:
Fix in advance a map $\psi$ from $\omega_1 \backslash \{0\} $
onto $\omega_1 \times \omega_1$,
such that $\alpha < \beta$ whenever $\psi(\beta) = (\alpha, \xi)$.
Now, given $X_p$, use CH and let
$\{F^p_\xi : \xi < \omega_1\}$ be a listing of all uncountable closed
subsets of $X_p$.
Whenever $0 < \beta < \omega_1$ and $\psi(\beta) = (\alpha, \xi)$
and $\lh(q) = \beta$, we may set $p = q \res \alpha$
and $F = F^p_\xi \subseteq X_p$.  It is sufficient to 
show how to accomplish (R17) with these specific $\alpha,\beta,p,q,F$.

Choose a perfect $K \subset F$ which is disjoint from
$\{\pi_p^{(q \res \zeta)}(u_{q \res(\zeta+1)}) :
\alpha \le \zeta < \beta \}$.
Then $\pi^q_p$ is 1-1 on $(\pi^q_p)\iv(K)$,
so choosing all  $u_{q\cat i}$ and 
$v^n_{q\cat i}$ in $(\pi^q_p)\iv (K)$ will ensure (R17).
Now \emph{fix} $i \in 2$, and write
$u$ and $v^n$  for
$u_{q\cat i}$ and $v^n_{q\cat i}$.
To ensure (R15) and (R9), we modify the argument of \cite{HK2}.
Let $\{Q^n : n \in \omega\}$ list $\FF_q$.
Let $d$ be a metric on $(\pi^q_p)\iv (K)$.
For each $s \in 2^{<\omega}$, choose a perfect
$L_s \subseteq (\pi^q_p)\iv (K)$.
Make these into a tree, in the sense that
each  $L_{s\cat 0} \cap L_{s\cat 1} = \emptyset$,
each $\diam(L_s) \le 2^{-\lh(s)}$,
and  $L_{s\cat 0} \subseteq L_s$ and $L_{s\cat 1} \subseteq L_s$.
Also make sure that whenever $\lh(s) = n+1$ we have
either $L_s \subseteq Q^n$ or $L_s \cap Q^n = \emptyset$.
Let $v[s\cat \ell]$ be any point in $L_{s \cat \ell} \backslash L_s$.
For $f \in 2^\omega$, let $\{u[f]\} = \bigcap_n L_{f\res n}$.
For any $f \in 2^\omega$, if we set $u = u[f]$
and $v^n = v[f\rest(n+1)]$, then (R15) will hold.
Now, (R9) requires $u \in B_{q\cat i}$.
Since $B_{q\cat i}$ is a Bernstein set and
$\{u[f] : f \in 2^\omega\}$ is a Cantor set, we may
choose $f$ so that  $u[f] \in B_{q\cat i}$.
\end{proofof}

If $H \subseteq X_g$ is closed and
for some  initial segment $p = g\res\alpha$ 
the projection $\pi^g_p(H) \in \FF_p$, 
then, by irreducibility,
$H = (\pi^g_p)\iv (\pi^g_p(H))$, so that $H$ is a $G_\delta$.
To make $X_g$ hereditarily Lindel\"of, it suffices to
capture projections for each closed $H \subseteq X_g$ this way,
but it is not clear
whether this can be done without using $\diamondsuit$.

\section{Remarks and Examples}
\label{sec-rem}
One cannot replace ``CH'' by
``$2^{\aleph_0} < 2^{\aleph_1}$'' in the statement of Theorem
\ref{thm-CH-weird}, since by Proposition \ref{prop-cov},
it is consistent with any cardinal arithmetic that
every non-scattered compactum of weight less than $\cccc$
contains a copy of the Cantor set.
As usual, define

\begin{definition}
$\cov(\MM)$ is the least $\kappa$ such that $\RRR$ is 
the union of $\kappa$ meager sets.
\end{definition}

Note that $\cov(\MM)$ is the least $\kappa$ such
that $\MA(\kappa)$ for countable partial orders fails.
Using this, we easily see:

\begin{lemma}
\label{lemma-nocov}
If $\kappa < \cov(\MM)$ and $E_\alpha \subset [0,1]$ is meager
for each $\alpha < \kappa$, then
$[0,1] \setminus \bigcup_{\alpha < \kappa} E_\alpha$ contains
a copy of the Cantor set.
\end{lemma}

\begin{proposition}
\label{prop-cov}
If $X$ is compact and not scattered,
and $w(X) < \cov(\MM)$, then $X$ contains a copy of the Cantor set.
\end{proposition}
\begin{proof}
Replacing $X$ by a subspace, we may assume that we have
an irreducible map $\pi: X \onto [0,1]$.
Let $\BB$ be an open base for $X$ with $|\BB| < \cov(\MM)$ and
$\emptyset \ne \BB$.

Whenever $U,V \in \BB$ with $\overline U \cap \overline V = \emptyset$,
let $E_{U,V} = \pi( \overline U) \cap \pi( \overline V)$.
Then $E_{U,V} \subset [0,1]$ is nowhere dense because
$\pi$ is irreducible.
Applying Lemma \ref{lemma-nocov}, let $K \subset [0,1]$ be a 
copy of the Cantor set disjoint from all the $E_{U,V}$.
Note that $|\pi\iv\{y\}| = 1$ for all $y \in K$.
Thus, $\pi\iv(K)$ is homeomorphic to $K$.
\end{proof}

Note that one can force ``$\cov(\MM) = \cccc$'' by adding
$\cccc$ Cohen reals, which does not change cardinal arithmetic,
but in the statement of Proposition \ref{prop-cov},
``$\cov(\MM)$'' cannot be replaced by ``$\cccc$''.
If CH holds in $V$, 
then one may force $\cccc$ to be arbitrarily large by
adding random reals,
and \textit{any} random real extension $V[G]$ 
will have a compact non-scattered
space of weight $\aleph_1$ which does not contain a Cantor subset.
In fact, Dow and Fremlin \cite{DF} show that 
if $X$ is a compact F-space in $V$, 
then in a random real extension $V[G]$, the corresponding compact space
$\widetilde X$ has no convergent $\omega$--sequences,
and hence no Cantor subsets.

The weird space constructed in \cite{HK2}
also failed to satisfy the CSWP 
(the complex version of the Stone--Weierstrass Theorem).  Using
the method there, we can modify the proof of Theorem \ref{thm-CH-weird} to get:

\begin{theorem}
\label{thm-CH-weirder}
Assuming CH, there is a separable first countable connected
weird space $X$ of weight $\aleph_1$ such that
$X$ fails the CSWP and $\KKK_X$ fails the CTP.
\end{theorem}
\begin{proof}
First, in the proof of Theorem \ref{thm-CH-weird},
replace $[0,1]$ by $\cdisc$, the closed unit disc in the complex plane,
so that we may view $X$
as a subspace of the $\aleph_1$--dimensional polydisc.
Then, as in \cite{HK2}, 
by carefully choosing the functions $h_{p\cat i}$, one can ensure
that the restriction to $X$ of 
the natural analog of the disc algebra refutes the CSWP
of $X$. 
To refute the CTP, construct in $\overline D$ a Cantor tree
$\{p_s : s \in 2^{<\omega}\} \subseteq \KKK_{\cdisc}$
such that  each $p_s$ is a wedge of the disc with center $0$
and radius $2^{-\lh(s)}$;  then each
$\bigcap_{n\in\omega}p_{f\res n} = \{0\}$.
Then, since we may assume  the point $0$ is not expanded 
in the construction of $X$,
the inverse images of the $p_s$ yield a counterexample to the CTP of $\KKK_X$.
\end{proof}

\end{document}